# Asymptotic formulae for the Lommel and Bessel functions and their derivatives

N.I. Aleksandrova

*N.A. Chinakal Institute of Mining, Siberian Branch, Russian Academy of Sciences,*

*Krasnyi pr. 54, Novosibirsk, 630091, Russia,*

*e-mail: alex@math.nsc.ru*

## 1. Summary

We derive new approximate representations of the Lommel functions in terms of the Scorer function and approximate representations of the first derivative of the Lommel functions in terms of the derivative of the Scorer function. Using the same method we obtain previously known approximate representations of the Nicholson type for Bessel functions and their first derivatives. We study also for what values of the parameters our representations have reasonable accuracy.

## 2. The main results

Solving problems in mechanics of discrete media [1, 2], we derived the following asymptotic formulae that are not available in the literature [3–14], but probably are of general interest:

$$s_{0,n}(ct) \approx -\frac{\pi}{2(ct/2)^{1/3}} \mathrm{Gi}\left[\frac{n-ct}{(ct/2)^{1/3}}\right], \quad (n=2k) \tag{2.1}$$

and

$$s'_{0,n}(ct) \approx \frac{\pi}{2(ct/2)^{2/3}} \mathrm{Gi}'\left[\frac{n-ct}{(ct/2)^{1/3}}\right], \quad (n=2k), \tag{2.2}$$

where $s_{0,n}$ is the Lommel function, Gi is the Scorer function, $c$ and $t$ are positive real numbers, $n$ and $k$ are positive integers. Prime denotes the derivative with respect to the argument. Each of formulae (2.1) and (2.2) holds true provided that $n \gg 1$ or $ct \gg 1$.

## 3. Motivation of our research

In [1, 2], the method of integral transformations was used to solve two-dimensional problems of wave propagation in discrete periodic media. In the process of solving those problems, it was necessary to find the functions $u$, $v$, provided that their Laplace–Fourier transforms are given by the formulae

$$u^{\mathrm{LF}}(p,q) = \frac{p}{p^2 + 4c^2 \sin^2(q/2)} \tag{3.1}$$

and

$$v^{LF}(p,q) = \frac{\sin(q/2)}{p^2 + 4c^2 \sin^2(q/2)},$$ (3.2)

where $c$ is the velocity of the propagation of disturbances, the superscript L denotes the Laplace transform (of parameter $p$) with respect to time $t$, and the superscript F denotes the discrete Fourier transform (of parameter $q$) with respect to $k$

$$f^L(p) = \int_0^\infty f(t)e^{-pt}dt \quad \text{and} \quad g^F(q) = \sum_{k=-\infty}^{k=\infty} g_k e^{iqk}.$$

Formally, the solution to the problem can be written as follows:

$$v_k(t) = \frac{1}{4\pi^2 i} \int_{-\pi}^{\pi} \int_{i\alpha-\infty}^{i\alpha+\infty} v^{LF}(p,q)e^{pt-iqk} dp dq.$$ (3.3)

A similar formula holds true for $u_k(t)$.

Inverting the Laplace transform [15], we obtain the following solutions:

$$u^F(t) = \cos\left[2ct \sin\left(\frac{q}{2}\right)\right]$$ (3.4)

and

$$v^F(t) = \frac{1}{2c} \sin\left[2ct \sin\left(\frac{q}{2}\right)\right].$$ (3.5)

Inverting the discrete Fourier transform in formulae (3.4), (3.5), we get:

$$u_k(t) = \frac{2}{\pi} \int_0^{\pi/2} \cos(2zk)\cos(2ct \sin z)dz = J_{2k}(2ct)$$ (3.6)

and

$$v_k(t) = \frac{1}{\pi c} \int_0^{\pi/2} \cos(2zk)\sin(2ct \sin z)dz = \frac{1}{\pi c} s_{0,2k}(2ct),$$ (3.7)

where $J_{2k}$ is the Bessel function of the first kind.

In the problems of mechanics [1, 2], it is important to be able to evaluate the behaviour of perturbations in the vicinity of the quasi-front $k = ct$ (quasi-front is a zone, where perturbations change from zero to maximum). Being motivated by this problem, we look for asymptotic representations of the Bessel and Lommel functions for $k \gg 1$.

In order to evaluate the behavior of function (3.6), we use the following asymptotic representation of the Bessel function:

$$J_n(ct) \approx \frac{1}{(ct/2)^{1/3}} \text{Ai}\left[\frac{n-ct}{(ct/2)^{1/3}}\right].$$ (3.8)

This formula is valid for $n \gg 1$ and is known as the Nicholson-type formula (see p. 142 in [9] or pp. 190 and 249 in [16]). Here,

$$\text{Ai}(z) = \frac{1}{\pi}\int_0^\infty \cos\left(zy + \frac{y^3}{3}\right)dy$$

is the Airy function.

We define

$$z = \frac{n-ct}{(ct/2)^{1/3}}. \tag{3.9}$$

Observe that it follows from (3.8) that the amplitude of $J_n(ct)$ in the neighbourhood of the point $ct = n$ (according to (11), this point can also be written as $z=0$) decreases as $t^{-1/3}$ (or $n^{-1/3}$) as $t \to \infty$ (or $n \to \infty$). Note also that the size of the zone, where $J_n(ct)$ increases from zero to the first maximum, increases as $t^{1/3}$ (or $n^{1/3}$).

Return to our mechanical problem. Substituting (3.8) into (3.6), we obtain the desired asymptotic representation for the function $u_k(t)$

$$u_k(t) \approx \frac{1}{(ct)^{1/3}} \text{Ai}\left[\frac{2(k-ct)}{(ct)^{1/3}}\right]. \tag{3.10}$$

Below we derive formulae (2.1) and (2.2) and study the limits of applicability of formulae (2.1), (2.2) and (3.8).

## 4. Derivation of formula (2.1)

Using the Slepyan method [16] of combined asymptotic ($t \to \infty$) inversion of the integral Laplace–Fourier transforms of long-wave disturbances in the vicinity of the ray $x = ct$, we can find the asymptotic behaviour for $v_k(t)$ that is similar to (3.10).

Applying the Slepyan method, we make the substitution $p = s + iq(c+c')$ and $k = (c+c')t$, where $c' \to 0$ and defines the vicinity of the ray $k = ct$, in the inner integral (3.3). This yields

$$v_k(t) = \frac{1}{4\pi^2 i}\int_{-\pi}^{\pi}\int_{i\alpha-\infty}^{i\alpha+\infty} v^{\text{LF}}(s+iq(c+c'),q)e^{st}dpdq.$$

We expand the numerator and denominator of the function $v^{\text{LF}}(s+iq(c+c'),q)$ in the Taylor series in a small neighbourhood of the point $q=0$ as $s \to 0$ and $c' \to 0$

$$v_k(t) \approx \frac{1}{4\pi^2 i}\int_{-\varepsilon}^{\varepsilon}\int_{\alpha-i\infty}^{\alpha+i\infty} \frac{\text{sign}(q)}{2ic(s+iqc'+iq^3c/24)}e^{st}dsdq,$$

where $\varepsilon > 0$ is small enough. Successively integrating and taking into account, that $c' = (k-ct)/t$, we obtain the following asymptotic formula that is similar to (12):

$$v_k(t) \approx -\frac{1}{4\pi c}\int_0^\varepsilon \sin\left(qc't + \frac{q^3 ct}{24}\right)dq \approx -\frac{1}{2c(ct)^{1/3}}\text{Gi}\left[\frac{2(k-ct)}{(ct)^{1/3}}\right], \tag{4.1}$$

where

$$\text{Gi}(z) = \frac{1}{\pi}\int_0^\infty \sin\left(zy + \frac{y^3}{3}\right)dy$$

is the Scorer function.

Comparing (3.7) and (4.1), we get the following approximate representation of the Lommel function $s_{0,n}$ for $n \gg 1$ in terms of the Scorer function Gi that is similar to (3.8):

$$s_{0,n}(ct) \approx -\frac{\pi}{2(ct/2)^{1/3}}\text{Gi}\left[\frac{n-ct}{(ct/2)^{1/3}}\right]. \tag{4.2}$$

Observe that the Lommel function $s_{0,n}$ is defined for even values of $n$ only (i.e. for $n = 2k$). It follows from (4.2) that the amplitude of $s_{0,n}(ct)$ in the neighbourhood of the point $ct = n$ decreases as $t^{-1/3}$ (or $n^{-1/3}$) as $t \to \infty$ (or $n \to \infty$). Note also the size of the zone, where $s_{0,n}(ct)$ decreases from zero to the first minimum, increases as $t^{1/3}$ (or $n^{1/3}$).

Finally, note that above we derived formula (3.10) from formula (3.8) of the Nicholson type solely for the sake of brevity. In fact, (3.10) can be obtained by using the Slepyan method of combined asymptotic inversion of the integral Laplace–Fourier transforms, just as we got above formula (4.1).

The approximate representation (4.2) of the Lommel function $s_{0,n}$ for $n \gg 1$ in terms of the Scorer function Gi is similar to the following formula (11.11.17) in [5]:

$$\mathbf{A}_{-\nu}(\nu + a\nu^{1/3}) \approx \frac{2^{1/3}}{\nu^{1/3}}\text{Hi}\left(-2^{1/3}a\right) + O(\nu^{-1}),$$

which gives an asymptotic expansion of the associated Anger–Weber function $\mathbf{A}_{-\nu}(z) = (1/\pi)\int_0^\infty \exp(\nu y - z\sinh y)dy$ for $\nu \to +\infty$ in terms of the Scorer function $\text{Hi}(z) = (1/\pi)\int_0^\infty \exp(zy - y^3/3)dy$.

## 5. Derivation of formula (2.2)

Observe that, in the following formula, the term $(n - ct)/(3t)$ can be neglected in a neighbourhood of the point $n = ct$ as $t \to \infty$:

$$\frac{dz}{dt} = \left(\frac{2}{ct}\right)^{1/3}\left(-c - \frac{n-ct}{3t}\right) \approx -c\left(\frac{2}{ct}\right)^{1/3}. \tag{5.1}$$

Differentiating (4.2) with respect to time, we get

$$\frac{ds_{0,n}(ct)}{dt} \approx -\frac{\pi}{2}\frac{d}{dt}\left[\frac{\text{Gi}(z)}{(ct/2)^{1/3}}\right] = -\frac{\pi}{2}\left(\frac{2}{ct}\right)^{1/3}\left[\frac{d\text{Gi}(z)}{dz}\frac{dz}{dt} - \frac{\text{Gi}(z)}{3t}\right].$$

Using (5.1) and assuming $t \to \infty$, we obtain the following asymptotic representation for the first derivative $s'_{0,n}$ for $n \gg 1$:

$$s'_{0,n}(ct) \approx \frac{\pi \text{Gi}'(z)}{2(ct/2)^{2/3}}. \tag{5.2}$$

Similarly, we derive the following asymptotic representation for the first derivative $J'_n$ for $n \gg 1$:

$$J'_n(ct) \approx -\frac{\text{Ai}'(z)}{(ct/2)^{2/3}}. \tag{5.3}$$

From (5.2) and (5.3), we conclude that, in a neighbourhood of the point $n = ct$, the functions $J'_n(ct)$ and $s'_{0,n}(ct)$ decrease as $t^{-2/3}$ (or $n^{-2/3}$) when $t$ (or $n$) increases. Note also the size of the zone, where $J'_n(ct)$ and $s'_{0,n}(ct)$ varies from zero to the first extremum, increases as $t^{1/3}$ (or $n^{1/3}$).

## 6. Numerical experiments

In order to determine the accuracy of the asymptotic representations (2.1), (2.2), (3.8) and (5.3) we plot the graphs of the functions appeared in (2.1), (2.2), (3.8) and (5.3). All figures given below are plotted for the case $c = 1$.

Let us use the following notations for the right-hand sides of formulae (2.1), (2.2), (3.8) and (5.3):

$$F_1(n,t) = \frac{1}{(t/2)^{1/3}} \text{Ai}\left[\frac{n-t}{(t/2)^{1/3}}\right], \qquad F_2(n,t) = -\frac{1}{(t/2)^{2/3}} \text{Ai}'\left[\frac{n-t}{(t/2)^{1/3}}\right],$$

$$F_3(n,t) = -\frac{\pi}{2(t/2)^{1/3}} \text{Gi}\left[\frac{n-t}{(t/2)^{1/3}}\right] \quad \text{and} \quad F_4(n,t) = \frac{\pi}{2(t/2)^{2/3}} \text{Gi}'\left[\frac{n-t}{(t/2)^{1/3}}\right].$$

In figures 1–4, we present the plots of the functions of the variable $t$, which appear in the left- and right-hand sides of formulae (2.1), (2.2), (3.8) and (5.3) for various values of $n$. The step of the variable $t$ is equal to 0.1. The dashed vertical lines correspond to the coordinates $t_* = n$ or $z = 0$.

Let us find the values of $n$, for which the accuracy of the asymptotic representations is reasonable. Let $\max|J_n(t)|$ denotes the maximum of the modulus of the function $J_n(t)$ calculated in a neighbourhood of the point circled in figure 1. The expressions $\max|J'_n(t)|$, $\max|s_{0,n}(t)|$,

$\max |s'_{0,n}(t)|$, $\max |F_1(n,t)|$, $\max |F_2(n,t)|$, $\max |F_3(n,t)|$ and $\max |F_4(n,t)|$ are defined similarly with the help of figures 2–4. Introduce the notation

$$\delta_1 = \left(1 - \frac{\max |F_1(n,t)|}{\max |J_n(t)|}\right) 100\%, \quad \delta_2 = \left(1 - \frac{\max |F_2(n,t)|}{\max |J'_n(t)|}\right) 100\%,$$

$$\delta_3 = \left(1 - \frac{\max |F_3(n,t)|}{\max |s_{0,n}(t)|}\right) 100\% \quad \text{and} \quad \delta_4 = \left(1 - \frac{\max |F_4(n,t)|}{\max |s'_{0,n}(t)|}\right) 100\%.$$

In tables 1 and 2, we give the values of the relative errors $\delta_1$, $\delta_2$, $\delta_3$ and $\delta_4$, calculated for the values of *n*, specified in figures 1–4. From tables 1 and 2, it follows that the relative errors $\delta_1$, $\delta_2$, $\delta_3$ and $\delta_4$ monotonically decrease as *n* increases.

**Table 1.** Relative errors $\delta_1$ and $\delta_2$.

| n | 2 | 6 | 10 | 20 |
|---|---|---|---|---|
| $\delta_1 (\%)$ | 4.7 | 2.7 | 2.1 | 1.4 |
| $\delta_2 (\%)$ | 7.4 | 5.2 | 4.1 | 2.7 |

**Table 2.** Relative errors $\delta_3$ and $\delta_4$.

| n | 6 | 10 | 20 | 40 |
|---|---|---|---|---|
| $\delta_3 (\%)$ | 9.7 | 7.4 | 5.0 | 3.4 |
| $\delta_4 (\%)$ | 6.2 | 5.1 | 3.6 | 1.8 |

From figures 1–4, for every pair of the functions, we see that, as *n* increases, the matching of the amplitudes of all local extrema get better, not only of the circled ones. From figures 1–4, we see also that, for every pair of the functions, the matching of the oscillation frequencies get better as *n* increases. For each pair of functions, the best approximation is achieved in a neighbourhood of the point $n = t$ (or $z = 0$). Note that, in the problems of mechanics, this neighbourhood corresponds to the quasi-front of the propagating wave.

Besides, we compared the plots of the functions (treated as functions of the variable *n*), which appear in the left- and right-hand sides of formulae (2.1), (2.2), (3.8) and (5.3), for various fixed values of *t*. That comparison showed that, for every pair of the functions, the plots agree best of all in a neighbourhood of the point $n = t$, that formulae (2.1), (2.2), (3.8) and (5.3) can be used starting from $t \approx 6$, and that the agreement gets better as *t* increases. For every pair of the functions, the distinction appears at a sufficiently large distance from the point $n = t$. Moreover, this distinction appears in the frequency of oscillations only. The maximal amplitudes of oscillations do not differ very much even at large distances from the point $n = t$.

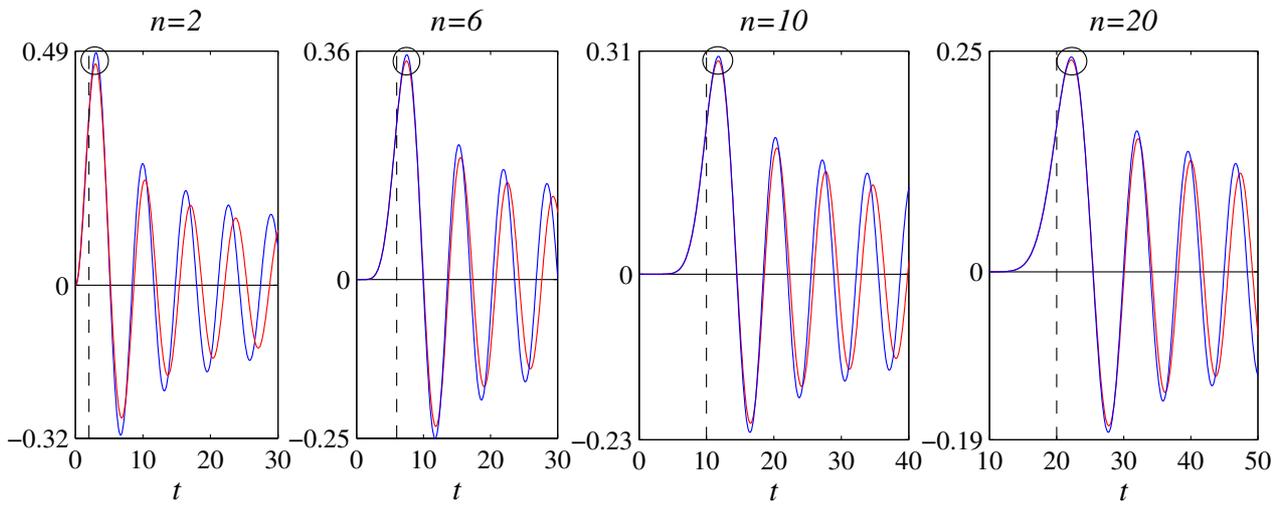

**Figure 1.** Plots of the functions: $J_n(t)$ – blue line; $F_1(n,t)$ – red line.

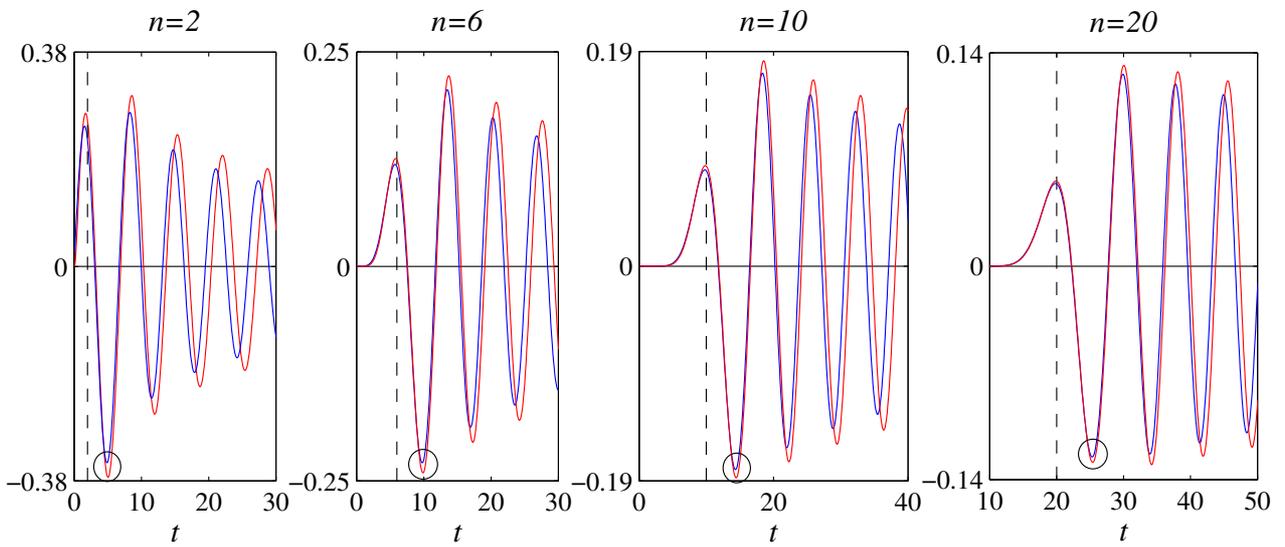

**Figure 2.** Plots of the functions: $J'_n(t)$ – blue line; $F_2(n,t)$ – red line.

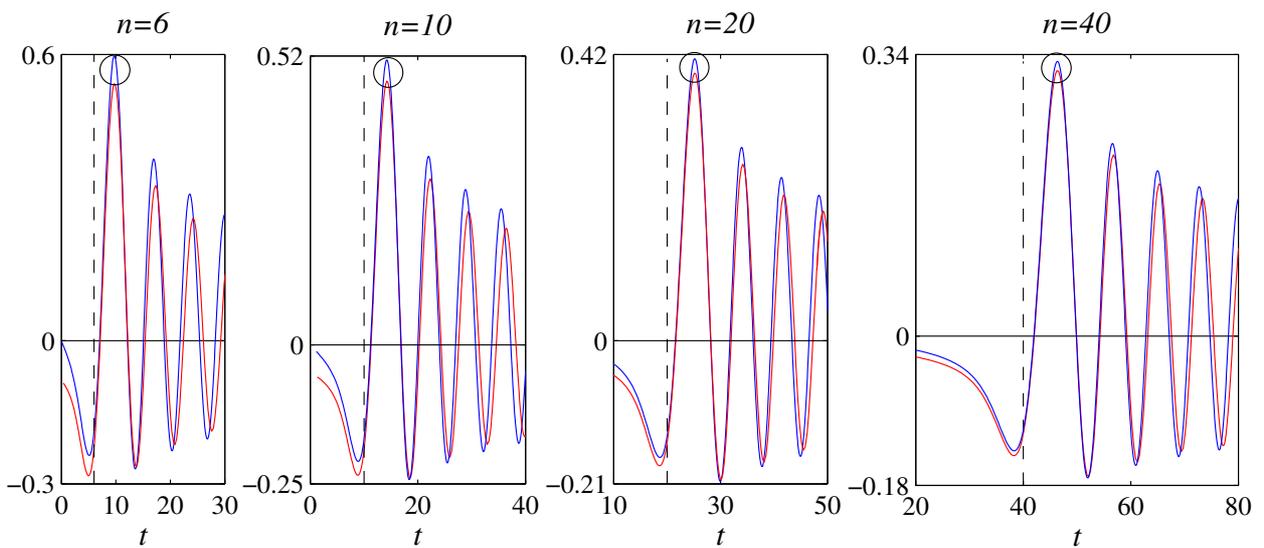

**Figure 3.** Plots of the functions: $s_{0,n}(t)$ – blue line; $F_3(n,t)$ – red line.

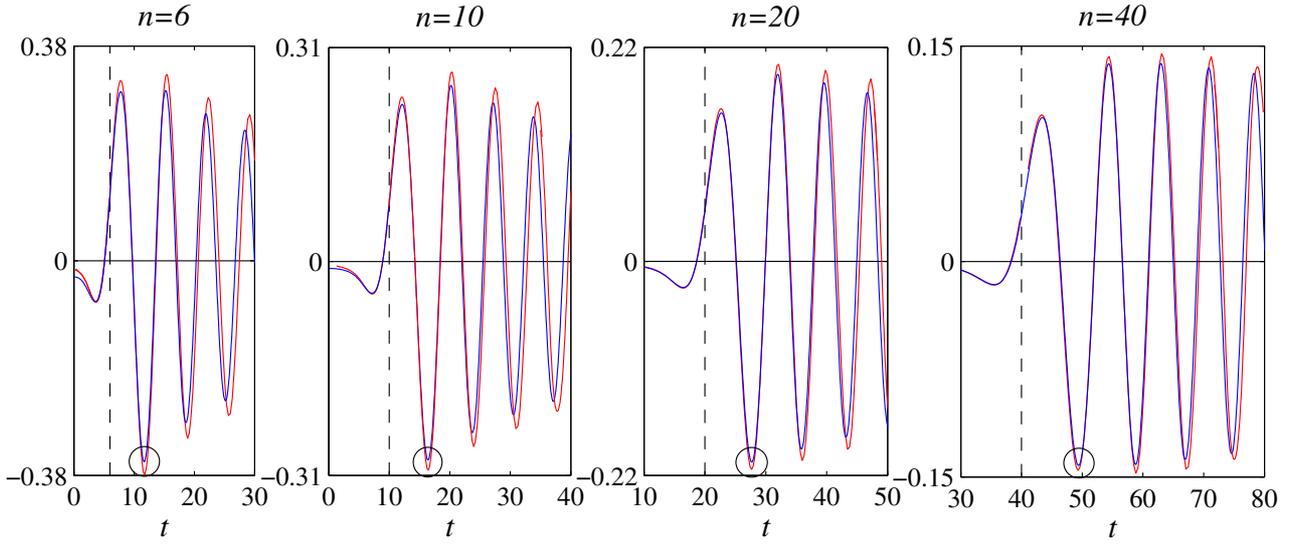

**Figure 4.** Plots of the functions: $s'_{0,n}(t)$ – blue line; $F_4(n,t)$ – red line.

## 7. Conclusion

As a result of the study of the asymptotic representations (2.1), (2.2), (3.8) and (5.3) it is shown that:

- the amplitudes of $s_{0,n}(ct)$ and $J_n(ct)$ in the neighbourhood of the point $ct = n$ decrease as $t^{-1/3}$ (or $n^{-1/3}$) as $t \to \infty$ (or $n \to \infty$);
- the amplitudes of $s'_{0,n}(ct)$ and $J'_n(ct)$ in the neighbourhood of the point $ct = n$ decrease as $t^{-2/3}$ (or $n^{-2/3}$) as $t \to \infty$ (or $n \to \infty$);
- the size of the zone, where $s_{0,n}(ct)$, $J_n(ct)$, $s'_{0,n}(ct)$ and $J'_n(ct)$ varies from zero to the first extremum, increases as $t^{1/3}$ (or $n^{1/3}$) as $t \to \infty$ (or $n \to \infty$) and
- representations (2.1), (2.2), (3.8) and (5.3) have reasonable accuracy starting from relatively small values of $n$ (namely, $n \approx 6$) or $t$ (namely, $ct \approx 6$).

## 8. Comparison of the results for the function $J_\nu(ct)$ described in [6] and this paper

Both [6] and this paper focus on the study of the behavior of the Bessel functions $J_\nu(ct)$ when the argument $ct$ and order $\nu$ are nearly equal. However, in this article, the emphasis is on the smallest values of $\nu = n$ and $ct$, for which the asymptotic formula (10) has reasonable accuracy for solving the problems of discrete periodic media [1, 2]. This differs our paper from [7], where the authors are looking for the values of $J_\nu(ct)$ for large values of $\nu$ and $ct$. In particular, Jentschura & Lötstedt [6] present, via apparently heroic numerical efforts, the following value $J_\nu(ct) = 0.002614463954691926$ for $\nu = 5000000.2$ and $ct = 5000000.1$. In this formula, the values of the argument $ct$ and order $\nu$ of the Bessel function are the largest ones for which we know the value of the Bessel function from the scientific literature. For $v = 5000000.2$ and $ct = 5000000.1$, the asymptotic formula (10) yields $J_v(ct) = 0.002614463961695188$. Hence, eight significant figures are in agreement with the exact numerical result given in [6].

In figure 5, we plot the graph of the function $F_1(v,t)$ for $v = 20000000.2$, which, according to formula (3.8), is asymptotically equivalent to the function $J_v(ct)$. From fig. 4 in [6] and figure 5, we see that the behaviour of the plots is the same.

Numerical experiments, discussed in this section, show that formula (3.8) is valid not only for the Bessel function $J_n(ct)$ of a positive integer order $n$, but also for the case where the order is a positive real number.

The problems associated with the confluence of the saddle points in the cusp region explained in [6] do not appear in our paper, because we use another method.

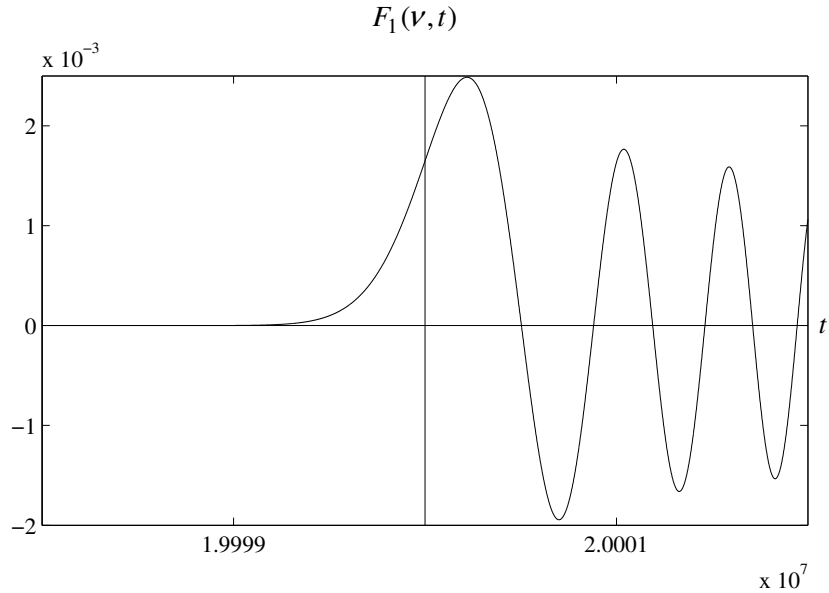

**Figure 5.** Plot of the function $F_1(v,t)$ for $v = 20000000.2$.

## 9. Comparison of the results for Bessel functions described in [9, 12] and this paper

Note that formulae (3.8) and (5.3) present the leading terms of the following complete asymptotic expansions of Bessel functions and their first derivatives (see pp. 281 and 287 in [12]):

$$J_v(y) \approx \sum_{k=0}^{\infty} (-1)^k \left(\frac{y}{2}\right)^{-(2k+1)/3} [P_k(\xi)\mathrm{Ai}(\xi) + Q_k(\xi)\mathrm{Ai}'(\xi)] \qquad (9.1)$$

and

$$J_v'(y) \approx \sum_{k=0}^{\infty} (-1)^{k+1} \left(\frac{y}{2}\right)^{-(2k+2)/3} [\overline{P}_k(\xi)\mathrm{Ai}(\xi) + \overline{Q}_k(\xi)\mathrm{Ai}'(\xi)]. \qquad (9.2)$$

Here $\xi = (\nu - y)(y/2)^{-1/3}$, $y \to \infty$, $|\arg y| \leq \pi - \varepsilon$, $y - \nu = O(y^{1/3})$, $\nu$ is a positive real number, the coefficients $P_k(\xi)$, $Q_k(\xi)$, $\overline{P}_k(\xi)$, $\overline{Q}_k(\xi)$ are polynomials in $\xi$ described in [12] (for example, $P_0(\xi) = 1$, $Q_0(\xi) = 0$, $\overline{P}_0(\xi) = 0$, $\overline{Q}_0(\xi) = 1$).

Besides formulae (3.8) and (5.3), (9.1) and (9.2), there are other asymptotic expansions of Bessel functions and their first derivatives for large values of the order. We mean the following formulae (10.19.8) and (10.19.12) in [5] (the same formulae are available on p. 414 in [9]):

$$J_\nu(\nu + a\nu^{1/3}) \approx \frac{2^{1/3}}{\nu^{1/3}} \mathrm{Ai}(-2^{1/3}a) \sum_{k=0}^{\infty} \frac{P_k(a)}{\nu^{2k/3}} + \frac{2^{2/3}}{\nu} \mathrm{Ai}'(-2^{1/3}a) \sum_{k=0}^{\infty} \frac{Q_k(a)}{\nu^{2k/3}} \qquad (9.3)$$

and

$$J'_\nu(\nu + a\nu^{1/3}) \approx -\frac{2^{2/3}}{\nu^{2/3}} \mathrm{Ai}'(-2^{1/3}a) \sum_{k=0}^{\infty} \frac{R_k(a)}{\nu^{2k/3}} + \frac{2^{1/3}}{\nu^{4/3}} \mathrm{Ai}(-2^{1/3}a) \sum_{k=0}^{\infty} \frac{S_k(a)}{\nu^{2k/3}}, \qquad (9.4)$$

where $a$ is a fixed complex number; $\nu$ is a complex number, such that $\nu \to \infty$ and its argument satisfies the inequality $|\arg \nu| \leq \pi/2 - \delta$ with some $\delta > 0$; the coefficients $P_k(a)$, $Q_k(a)$, $R_k(a)$, $S_k(a)$ are polynomials in $a$, in particular, $P_0(a) = 1$, $Q_0(a) = 3/10\,a^2$, $R_0(a) = 1$ and $S_0(a) = 3/5\,a^3 - 1/5$.

At first glance, formulae (9.1) and (9.2) are similar to formulae (9.3) and (9.4). But in fact, they differ radically, because in (9.3) and (9.4) the expansion is carried out in powers of the order $\nu$ of the functions, whereas in formulae (9.1) and (9.2) the expansion is carried out in the powers of the argument of the functions.

**Acknowledgement**

The author would like to thank the anonymous referees for their helpful comments.